\documentclass{gtart}

\def\ifplaintex{\expandafter\ifx\csname documentclass\endcsname\relax}

\def\gtp{{\mathsurround=0pt\it $\cal G\mskip-2mu$eometry \&\ 
$\cal T\!\!$opology $\cal P\!$ublications}}  

\def\recd{{\small Received:\qua\receiveddate\ifx\reviseddate\relax
\else\qquad Revised:\qua\reviseddate\fi\par}} 

\def\lognumber#1{\def\thelognumber{#1}}
\def\volumenumber#1{\def\thevolumenumber{#1}}
\def\volumeyear#1{\def\thevolumeyear{#1}}
\def\papernumber#1{\def\thepapernumber{#1}}
\def\pagenumbers#1#2{\def\startpage{#1}\def\finishpage{#2}}
\def\published#1{\def\publishdate{#1}}

\def\received#1{\def\receiveddate{#1}}
\def\revised#1{\def\reviseddate{#1}}
\def\accepted#1{\def\accepteddate{#1}}

\def\asciiaddress#1{\def\theasciiaddress{#1}}
\def\asciiemail#1{\def\theasciiemail{#1}}

\long\def\asciiabstract#1{\long\def\theasciiabstract{#1}}

\let\\\par\let\thelognumber\relax\let\thevolumenumber\relax
\let\thepapernumber\relax\let\thevolumeyear\relax\let\startpage\relax
\let\finishpage\relax\let\publishdate\relax\let\receiveddate\relax
\let\reviseddate\relax\let\accepteddate\relax\let\theasciititle\relax
\let\theasciiauthors\relax\let\theasciiaddress\relax
\let\theasciiabstract\relax

\let\theasciiemail\relax

\ifplaintex
\font\logobig=cmssbx10 scaled 3836
\font\logomed=cmssbx10 scaled 2557
\else
\font\logobig=cmssbx10 scaled 4200
\font\logomed=cmssbx10 scaled 2800
\fi

\long\def\makeagttitle{   
\count0=\startpage
\agt\hfill      
\hbox to 45truept{\vbox to 0pt{\vglue -13truept{\logomed A\kern -.37em{\logobig 
T}\kern -.38em G}\vss}\hss}
\break
{\small Volume \thevolumenumber\ (\thevolumeyear)
\startpage--\finishpage\nl
Published: \publishdate}

\vglue .25truein

{\parskip=0pt\leftskip 0pt plus
1fil\def\\{\par\smallskip}{\Large\bf\thetitle}\par\medskip} \vglue
0.05truein

{\parskip=0pt\leftskip 0pt plus 1fil\def\\{\par}{\sc\theauthors}
\par\medskip}%
 
\vglue 0.03truein 

{\small\leftskip 25truept\rightskip 25truept{\bf Abstract}\stdspace\theabstract

{\bf AMS Classification}\stdspace\theprimaryclass
\ifx\thesecondaryclass\relax\else; \thesecondaryclass\fi\par
{\bf Keywords}\stdspace \thekeywords\par}\vglue 7truept

}   

\ifplaintex
\font\phead=cmsl9 scaled 950
\font\pnum=cmbx10 scaled 913
\font\pfoot=cmsl9 scaled 950
\headline{\vbox to 0pt{\vskip -4.5mm\line{\small\phead\ifnum
\count0=\startpage ISSN 1472-2739 (on-line) 1472-2747 (printed)
\hfill {\pnum\folio}\else\ifodd\count0\def\\{ }%
\ifx\theshorttitle\relax\thetitle\else\theshorttitle\fi\hfill{\pnum\folio}
\else\def\\{ and }{\pnum\folio}\hfill\ifx\theshortauthors\relax\theauthors
\else\theshortauthors\fi\fi\fi}\vss}}
\footline{\vbox to 0pt{\vglue 0mm\line{\small\pfoot\ifnum\count0=\startpage
\copyright\ \gtp\hfill\else
\agt, Volume \thevolumenumber\ (\thevolumeyear)\hfill\fi}\vss}}
\else
\font=cmsl9 scaled 1050
\font\lnum=cmbx10 
\font=cmsl9 scaled 1050
\makeatletter
\def\@oddhead{{\small\lhead\ifnum\count0=\startpage ISSN 1472-2739 
(on-line) 1472-2747 (printed)\hfill {\lnum\number\count0}\else\ifodd\count0
\def\\{ }\ifx\theshorttitle\relax \thetitle \else\theshorttitle\fi\hfill
{\lnum\number\count0}\else\def\\{ and }{\lnum\number\count0}
\hfill\ifx\theshortauthors\relax 
\theauthors\else\theshortauthors\fi\fi\fi}}\def\@evenhead{\@oddhead}
\def\@oddfoot{\small\lfoot\ifnum\count0=\startpage\copyright\ \gtp\hfill\else
\agt, Volume \thevolumenumber\ (\thevolumeyear)\hfill\fi}
\def\@evenfoot{\@oddfoot}
\makeatother
\fi
\let\maketitlepage\makeagttitle
\let\makeshorttitle\maketitlepage
\let\maketitle\maketitlepage

\newwrite\gtoutfile
\long\gdef\makeheadfile{  
{\def\\{, }\def\s{ }
\immediate\openout\gtoutfile head.xxx
\immediate\write\gtoutfile{Proxy-for: \ifx\theasciiauthors\relax
\theauthors\else\theasciiauthors\fi\s<\ifx\theasciiemail\relax\theemail\else\theasciiemail\fi>}
\immediate\write\gtoutfile{\noexpand\\}
\immediate\write\gtoutfile{Authors: \ifx\theasciiauthors\relax
\theauthors\else\theasciiauthors\fi}
{\def\\{ }\immediate\write\gtoutfile{Title: \ifx\theasciititle\relax
\thetitle\else\theasciititle\fi}}
\immediate\write\gtoutfile{Subj-class: GT or SG, GR etc}
\immediate\write\gtoutfile{MSC-class: \theprimaryclass\ifx\thesecondaryclass\relax\else, \thesecondaryclass\fi}
\immediate\write\gtoutfile{Journal-ref: Algebr. Geom. Topol. \thevolumenumber\s
(\thevolumeyear) \startpage-\finishpage}
\immediate\write\gtoutfile{Comments: Published by Algebraic and
Geometric Topology at}
\immediate\write\gtoutfile{\s\s\s  http://www.maths.warwick.ac.uk/agt/AGTVol\thevolumenumber/agt-\thevolumenumber-\thepapernumber.abs.html}
\immediate\write\gtoutfile{\noexpand\\}
\immediate\write\gtoutfile{}
\ifx\theasciiabstract\relax
\immediate\write\gtoutfile{\theabstract}\else
\immediate\write\gtoutfile{\theasciiabstract}\fi
\immediate\write\gtoutfile{}
\immediate\write\gtoutfile{\noexpand\\}
\immediate\write\gtoutfile{}
\immediate\closeout\gtoutfile}}  

\def\maketitlepage{\makeagttitle\makeheadfile}
\let\makeshorttitle\maketitlepage
\let\maketitle\maketitlepage

\lognumber{30}
\volumenumber{3}
\volumeyear{2003}
\papernumber{30}
\published{29 September 2003}
\pagenumbers{905}{920}
\received{23 October 2002}
\revised{5 July 2003}
\accepted{5 September 2003}

\usepackage{amsmath, amssymb, epsfig}

\newtheorem{theorem}{Theorem}[section]    
\newtheorem{lemma}[theorem]{Lemma}          
   
\newtheorem{proposition}[theorem]{Proposition}   
\theoremstyle{remark}   
\newtheorem{remark}[theorem]{Remark }   
\newtheorem{definition}[theorem]{Definition}

\newcommand{\Cc}{{C}}
\newcommand{\Nn}{\mathbb{N}}   
\newcommand{\Rr}{\mathbb{R}}

\newcommand{\Zz}{\mathbb{Z}}   
\newcommand{\Qq}{\mathbb{Q}}   
 
\DeclareMathOperator{\Sign}{Sign}
\DeclareMathOperator{\Hom}{Hom}

\begin{document} 
\title{On the slice genus of  links} 
\authors {Vincent Florens\\Patrick M. Gilmer} 
\address{Laboratoire I.R.M.A. Universit\'e Louis Pasteur\\Strasbourg, France }
\secondaddress{Department of Mathematics,  
  Louisiana State University\\  
 Baton Rouge, LA 70803, USA}  
\asciiaddress{Laboratoire I.R.M.A. Universite Louis Pasteur\\Strasbourg, France\\and\\Department of Mathematics,  
  Louisiana State University\\  
 Baton Rouge, LA 70803, USA}  

\email{vincent.florens@irma.u-strasbg.fr}
\secondemail{gilmer@math.lsu.edu} 
\asciiemail{vincent.florens@irma.u-strasbg.fr, gilmer@math.lsu.edu} 

\begin{abstract}
We define  Casson-Gordon $\sigma$-invariants for links and  
 give a lower bound of the slice genus of a link in terms of these 
invariants.  
 We study as an example a family
of two component links of genus $h$
 and show that their slice genus is $h$,  whereas  
 the Murasugi-Tristram  inequality does not obstruct this link from bounding 
an 
annulus
 in the 4-ball.
 \end{abstract}

\asciiabstract{We define Casson-Gordon sigma-invariants for links and
give a lower bound of the slice genus of a link in terms of these
invariants.  We study as an example a family of two component links of
genus h and show that their slice genus is h, whereas the
Murasugi-Tristram inequality does not obstruct this link from bounding
an annulus in the 4-ball.}
 
\primaryclass{57M25}\secondaryclass{57M27}
\keywords {Casson-Gordon invariants, link signatures}

\maketitle  

\section{Introduction} 
  
A knot in $S^3$ is slice if it bounds a smooth $2$-disk
 in the $4$-ball $B^4$. 
Levine showed \cite{Le} that a slice knot is 
algebraically
 slice, i.e.\ any Seifert form of
a slice
 knot is metabolic. In this case, 
the Tristram-Levine signatures 
at
 the
prime power order roots of 
unity 
of 
a slice
knot must be zero.
Levine showed also that the 
converse holds
 in high odd dimensions, i.e.\
 any algebraically slice knot is slice.
This is false in dimension $3$:  
 Casson and Gordon \cite{CG1, CG2, G} showed that certain   
two-bridge 
 knots in 
$S^3$, which are algebraically slice,
 are not slice knots.
 For this purpose, they defined several knot and 3-manifold invariants, closely related to the Tristram-Levine 
 signatures of associated links.
 Further methods to calculate these invariants were developed by Gilmer \cite{Gi3,Gi4},   
 Litherland \cite{Li}, Gilmer-Livingston \cite{GL}, and Naik \cite{N}. Lines \cite{L} also computed some of these invariants for some fibered knots, which are algebraically slice
  but not slice.
The slice genus of a link is the minimal genus 
for a smooth
oriented connected surface
properly embedded
 in $B^4$
 with boundary the given link.

The Murasugi-Tristram inequality (see Theorem \ref{MT} below) gives a lower bound on the slice genus of a link in terms of  
the link's
Tristram-Levine signatures
and
related nullity invariants.
 The second author
 \cite{Gi1} used Casson-Gordon invariants to give another lower bound on the slice genus of a knot. In particular he gave examples of algebraically slice knots whose slice genus is arbitrarily large.  We apply these methods to restrict the slice genus of a link.  
  
We study as an example a family of two component links, which have genus $h$  Seifert surfaces.  Using Theorem \ref{main}, we show that these links cannot bound a smoothly embedded 
surface
 in $B^4$ with genus lower than $h$, while the Murasugi-Tristram 
inequality does not show this. In fact there are some links with the same Seifert 
form
that bound 
annuli
 in  
 $B^4$. We work in the smooth category.

 The second author was partially supported by  NSF-DMS-0203486.
   
\section{Preliminaries}   
   
\subsection{The Tristram-Levine signatures}

Let $L$ be an oriented link in $S^3$, with $\mu$ 
components,
 and  
 $\theta_S$ be  
 the Seifert pairing corresponding to  
a connected Seifert surface $S$ of the link.  
For any complex number $\lambda$ with $| \lambda | =1$, one considers the hermitian form  
$\theta_S^\lambda := (1-\lambda) \theta_S + (1 - \overline{\lambda})  
 (\theta_S)^T$.
The Tristram  signature $\sigma_L(\lambda)$ and nullity $n_L(\lambda)$  
 of $L$ are defined as the signature and nullity of $\theta_S^\lambda$.  
Levine 
defined these same
 signatures for knots \cite{Le}.  
The Alexander polynomial of $L$ is 
 $\Delta_L (t) := \text{Det} ( \theta_S - t (\theta_S)^T).$  
As is well-known, $\sigma_L$ is
a locally constant map on the complement in $S^1$  
of the roots of $\Delta_L$ and $n_L$ is zero on this complement.  
If  $\Delta_L = 0, $ it is still true that the signature and nullity are locally constant functions  
on the complement of some finite collection of points.  
  
The Murasugi-Tristram  inequality allows one to estimate 
the slice genus of $L$, in terms of the values of $\sigma_L(\lambda)$ and $n_L(\lambda)$.
\begin{theorem}{\rm \cite{M,T}}\qua   
Suppose that $L$ is the boundary of a 
properly embedded
connected oriented surface $F$ of genus $g$
 in $B^4$. Then, 
 if $\lambda$ is a prime power order 
root of 
unity,  
we have   
$$ |\sigma_L(\lambda)| + n_L(\lambda) \leq 2 g + \mu - 1. $$
\label{MT}   
\end{theorem}

\subsection{The Casson-Gordon $\sigma$-invariant}

In this section, for the reader 
convenience, we review the definition  and some of the properties of the simplest  
kind of Casson-Gordon 
invariant.
 It is a reformulation of the Atiyah-Singer $\alpha$-invariant.    

Let $M$ be an oriented compact 
 three manifold and $\chi \co H_1(M) \rightarrow\mathbb {C}^*$   
be a character of finite order.
For some $q \in \Nn^*$, the image of $\chi$ 
is contained  a cyclic subgroup of order $q$ 
generated by $\alpha=e^{2i \pi /q}$. 
As $\Hom(H_1(M), \Cc_q) =[M,B(\Cc_q)]$, 
 it follows that   
  $\chi$ induces $q$-fold covering of $M$, denoted $\widetilde{M}$,  
with a canonical deck transformation. 
We will denote this transformation also by $\alpha.$
If $\chi$ maps onto  $\Cc_q,$ 
the canonical deck transformation sends $x$ to the other endpoint of the arc
  
that begins at $x$   
 and covers  a loop representing an element of
  $(\chi)^{-1}(\alpha)$.    
   
As the bordism group $\Omega_3(B(\Cc_q))= \Cc_q$, 
 we may conclude that  
 $n$ disjoint copies of $M$ , for some integer $n$,  bounds
 bound 
 a compact $4$-manifold $W$ over $B(\Cc_q)$.  
  Note $n$ can be taken
to be $q.$
 Let $\widetilde{W}$ be the induced covering with the deck transformation, denoted also by $\alpha$, that restricts   
 to $\alpha$ on the boundary.   
 This
 induces a $\Zz[\Cc_q]$-   
module structure on  $C_*(\widetilde{W})$, 
where the multiplication by $\alpha \in \Zz[\Cc_q]$   
 corresponds to the action of $\alpha$ on $\widetilde{W}.$

The cyclotomic field $\Qq(\Cc_q)$ is a natural $\Zz[\Cc_q]$-module and  
the twisted homology $H_*^t(W;\Qq(\Cc_q))$ is defined as the homology   
 of  $$C_*(\widetilde{W}) \otimes_{\Zz[\Cc_q]} \Qq(\Cc_q). $$  
Since $\Qq(\Cc_q)$ is   
 flat over $\Zz[\Cc_q]$, we get an isomorphism $$ H_*^t (W;\Qq(\Cc_q)) \simeq H_*(\widetilde{W})    
 \otimes_{\Zz[\Cc_q]} \Qq(\Cc_q). $$ Similarly,   
the twisted homology $H_*^t(M;\Qq(\Cc_q))$ is defined as the homology   
 of  $$C_*(\widetilde{M}) \otimes_{\Zz[\Cc_q]} \Qq(\Cc_q). $$   
Let $\widetilde{\phi}$ be the intersection form on $H_2(\widetilde{W};\Qq)$   
 and  define $$\phi_{\chi} (W)\co   
  H_2^t(W;\Qq(\Cc_q)) \times H_2^t(W;\Qq(\Cc_q))   
  \to \Qq(\Cc_q) $$   
so that, 
for all $a,b$ in $\Qq(\Cc_q)$ 
and $x,y$ in $H_2(\widetilde{W})$, 
 $$ \phi_{\chi} (W)(x \otimes a, y \otimes b) = \overline{a} b \sum_{i=1}^q   
\widetilde{\phi} (x,\alpha^i y) \overline{\alpha}^i, $$  
 where $a \rightarrow \bar{a}$ denotes the involution on
$\Qq(\Cc_q)$ induced by 
 complex conjugation.

\begin{definition}   
The Casson-Gordon $\sigma$-invariant of $(M,\chi)$ and the related nullity are
$$ \sigma(M,\chi):= \frac{1}{n} \big( {\Sign}(\phi_{\chi}(W)) - {\Sign}(W) \big) $$    
$$ \eta(M,\chi):= {\dim} \  H^t_1(M;\Qq(\Cc_q)). $$   
\label{sigma}  
\end{definition}

If $U$ is a closed 4-manifold and 
$\chi \co H_1(U) \rightarrow C_q$  we may define $\phi_{\chi}(U) $
as above. One has that modulo torsion the bordism group  
$\Omega_4(B(\Cc_q)) $ is generated by the constant map from $CP(2)$ to  $B(\Cc_q).$ 
If $\chi$ is trivial, one has that ${\Sign}(\phi_{\chi}(U))= {\Sign}(U). $ Since both signatures are invariant under cobordism, 
 one has  in  general that ${\Sign}(\phi_{\chi}(U))= {\Sign}(U).$  
The independence of  $\sigma(M,\chi)$ from the choice of $W$ 
 and $n$ follows from this and
 Novikov additivity. One may see directly that 
 these invariants do not depend on the choice of $q$. In this way Casson and Gordon argued that   $\sigma(M,\chi)$ is an invariant. Alternatively one may use the Atiyah-Singer G-Signature  
 theorem and 
Novikov additivity  \cite{AS}.

We now describe a way
to compute  $ \sigma(M,\chi)$ for a given
 surgery presentation of $(M,\chi)$.  
  
\begin{definition}   
Let $K$ be an oriented knot in $S^3$. Let $A$ be an embedded annulus such that $\partial A=K   
 \cup K^{\prime}$ with $lk(K,K^{\prime})=f$. 
A 
 \emph{p-cable on $K$ with twist $f$} is   
 defined 
to be
 the union of oriented parallel copies of $K$ lying in $A$ such that the    
 number of copies with the same orientation minus the number with opposite orientation is   
 equal to $p$.    
\label{cable}   
\end{definition}

Let us suppose that $M$ is obtained by surgery on a framed link $L=L_1 \cup \dots \cup L_\mu$   
 with framings $f_1,\dots,f_\mu$. One shows that   
 the linking matrix $\Lambda$ of $L$ with framings in the diagonal is a presentation matrix   
 of $H_1(M)$ and a character on $H_1(M)$ is determined by 
${\alpha}^{p_i}=\chi(m_{L_i}) \in \Cc_q$   
 where $m_{L_i}$  
denotes  
 the class of the meridian of $L_i$.   
Let $\vec p= (p_1, {\dots}, p_\mu)$. 
 We use the following generalization 
of a formula in \cite[Lemma (3.1)]{CG2}, 
where all $p_i$ are assumed to be 
 $1$, that is  given in \cite[Theorem(3.6)]{Gi2}.

\begin{proposition} \label {surgeryformula}  
Suppose $\chi$ maps onto $C_q$.
Let $L^{\prime}$ with $\mu^\prime$ components  
 be the link obtained from $L$ by replacing each component by a non-empty     
algebraic $p_i$-cable with twist $f_i$ along this component.   
Then, if $\lambda=e^{2i r \pi /q}$, for $(r,q)=1$, one has   
$$ \sigma(M,\chi^r) = \sigma_{L^{\prime}}(\lambda) - \Sign(\Lambda) +    
2  \frac{r(q-r)}{q^2} \vec p^\top \Lambda \vec p, $$   
$$ \eta(M,\chi^r) = \eta_{L^{\prime}}(\lambda)-\mu^\prime + \mu. $$  
\end{proposition}

The following proposition collects
some easy
additivity properties of
 the $\sigma$-invariant and the nullity under the connected sum.

\begin{proposition} \label {add} 
Suppose that $M_1,$ $M_2$ are connected. 
Then, \\ for all  $\chi_i \in H^1(M_i; C_q)$, $i=1,2$,  we have
$$\sigma(M_1\# M_2, \chi_1 \oplus \chi_2) = \sigma(M_1,\chi_1) +  \sigma(M_2,\chi_2). $$
If both $\chi_i$ are
non-trivial, then
$$\eta(M_1\# M_2,\chi_1 \oplus \chi_2) = \eta(M_1,\chi_1) +  \eta(M_2,\chi_2) + 1.$$
If one $\chi_i$ is trivial,
then
$$\eta(M_1\# M_2,\chi_1 \oplus \chi_2) = \eta(M_1,\chi_1)+  \eta(M_2,\chi_2). $$
\end{proposition}

\begin{proposition} \label {s1s2} 
For all $\chi \in H_1(S^1 \times S^2; C_q)$,
we have
$$ \sigma(S^1 \times S^2, \chi) = 0 $$ 
$$ \text{If } \chi \neq 0 \text{, then } \eta(S^1 \times S^2,\chi)=0. \text{ If } \chi=0 \text{, then } \eta(S^1 \times S^2,\chi)= 1 .$$  
\end{proposition}

 Proposition \ref{s1s2}  for non-trivial 
 $\chi$ can be proved for example
 by the use of Proposition \ref{surgeryformula}, since $S^1 \times S^2$ is
 obtained by surgery on the unknot framed $0$. However it is simplest to derive this result directly from the 
definitions.

\subsection{The Casson-Gordon $\tau$-invariant}   \label{cgt}

In this section, we recall the definition and some of the properties of the Casson-Gordon 
$\tau$-invariant.
Let  $C_\infty$ 
denote
 a multiplicative infinite cyclic group generated by $t.$   
  For $\chi^+  \co H_1(M) \to \Cc_q \oplus \ C_\infty$, we denote $\bar{\chi}  \co H_1(M) \to \Cc_q $ the character obtained by composing $\chi^+$ with projection on the first factor.  The character $\chi^+$ induces a $\Cc_q \times  C_\infty$-covering $\widetilde{M}_\infty$ of $M$.   
   
Since the bordism group $\Omega_3(B(\Cc_q \times  C_\infty ))= \Cc_q,$ 
 bounds a compact $4$-manifold $W$ over $B(\Cc_q \times  C_\infty)$ 
  Again $n$ can be  taken from
to be $q$.
   
 If we identify $\Zz[\Cc_q \times  C_\infty]$ with the    
 Laurent polynomial ring $\Zz[\Cc_q][t,t^{-1}]$, the field $\Qq(\Cc_q)(t)$ of rational functions   
 over the cyclotomic field $\Qq(\Cc_q)$ is a flat  $\Zz[\Cc_q \times C_\infty]$-module.   
We consider the chain complex $C_*(\widetilde{W}_\infty)$ as a  $\Zz[\Cc_q \times C_\infty]$-module   
 given by the deck transformation of the covering.   
 Since $W$ is compact, the vector space $H_2^t(W;\Qq(\Cc_q)(t)) \simeq H_2(\widetilde{W}_\infty)   
 \otimes_{\Zz[\Cc_q][t,t^{-1}]} \Qq(\Cc_q)(t)$ is finite dimensional.

We let $J$ denote the  involution 
 on  
$\Qq(\Cc_q)(t)$ that is linear over $\Qq$
sends $t^i$  to $t^{-i}$ and 
 $\alpha^i$  to $\alpha^{-i}.$
As in   \cite{G},  
 one defines
a hermitian form, with respect to $J$,
$$\phi_{\chi^+}  \co H_2^t(W;\Qq(\Cc_q)(t)) \times   
  H_2^t(W;\Qq(\Cc_q)(t)) \to \Qq(\Cc_q)(t), $$ such that
 $$ \phi_{\chi^+} (x \otimes a, y \otimes b) = J(a) \cdot b \cdot \sum_{i \in \Zz}   \sum_{j= 1} ^q  \widetilde{\phi^+} (x,t^i \alpha^j  y)
  \overline{\alpha}^j t^{-i}.$$
Here $\widetilde{\phi^+} $  
 denotes  
 the ordinary intersection form on $\widetilde{W}_\infty.$
Let $\mathcal{W} (\Qq(\Cc_q)(t))$ be the Witt group of non-singular hermitian forms on finite   
 dimensional $\Qq(\Cc_q)(t)$ vector spaces.  
 Let us consider $H_2^t(W;\Qq(\Cc_q)(t)) / (\text{Radical}( \phi_{\chi^+}))$.  
  The induced form   
 on it represents an element in $\mathcal {W}$ $ (\Qq(\Cc_q)(t)),$ which we denote 
$w(W)$.
 Furthermore, the ordinary intersection form on $H_2(W;\Qq)$ represents an element of $\mathcal{W} (\Qq)$.
 Let $w_0(W)$ be the image of this element in  $\mathcal{W} (\Qq(\Cc_q)(t))$.

\begin{definition}   
 The Casson-Gordon $\tau$-invariant of $(M,\chi^+)$ is 
 $$ \tau (M,\chi^+)  : = \frac{1}{n} \big( w(W) - w_0(W) \big) 
\in \mathcal{W}   (\Qq(\Cc_q)(t)) \otimes \Qq. $$   
\end{definition}

Suppose that $nM$ bounds another
 compact $4$-manifold $W^\prime$ over $B(\Cc_q \times  C_\infty)$. Form the closed
 compact manifold over $B(\Cc_q \times  C_\infty)$, $U := W \cup W^\prime$ by gluing along
 the boundary. 
By Novikov additivity, 
we get $w(U)-w_0(U)= \big( w(W) - w_0(W) \big) - \big( w(W^\prime) - w_0(W^\prime) \big)$.
Using \cite{CF},  
the bordism group $\Omega_4(B(\Cc_q \times  C_\infty)) $, modulo torsion, 
is generated by $CP(2)$, with the constant map to  $B(\Cc_q \times  C_\infty)$.  
We have that
 $w( CP(2))=w_0(CP(2))$.
Since  $w(U)$, and $w_0(U)$ only depend on the bordism class of $U$ over $B(\Cc_q \times  C_\infty)$, it follows that $w(U)=w_0(U)$ and $ \tau (M,\chi^+)$ is independent
 of the choice of $W$.  Using the above techniques, one may check $\tau (M,\chi^+)$  is  independent of $n$.

If $A \in \mathcal{W}  (\Qq(\Cc_q)(t)),$ let $A(t)$ be a matrix representative for $A$. The entries of $A(t)$   
 are Laurent polynomials with coefficients in $\Qq(\Cc_q)$. If 
 $\lambda$ is in 
$S^1 \subset \mathbb{C}$,
then   
 $A(\lambda)$ is hermitian and has a well defined signature $\sigma_\lambda(A)$.    
 One can view $\sigma_\lambda(A)$ as a locally constant map
 on the complement   
 of 
the
set of
 the zeros of 
 $\det A(\lambda)$.
 As in \cite{CG1}, we re-define $\sigma_\lambda(A)$ at each point of discontinuity   
 as the average of the one-sided limits at the point.     

We have the following estimate \cite[Equation (3.1)]{Gi3}.

\begin{proposition}\label{firstest} Let $\chi^+ \co H_1(M) \to \Cc_q \oplus C_\infty$ and $\bar \chi \co H_1(M) \to \Cc_q$  be $\chi^+$ followed by the projection to $\Cc_q$.
 We have   
$$ | \sigma_{1} \big(\tau (M,\chi^+) \big) - \sigma(M,\bar {\chi}) | \leq \eta(M,\bar{\chi}).$$
\end{proposition}

\subsection{Linking forms}   
   
  Let $M$ be a rational homology 3-sphere with linking form $$l  \co  H_1(M) \times H_1(M) \to   
 \Qq/\Zz.$$   
We have that  $l$ is non-singular, 
that is  the adjoint of $l$  is  an isomorphism $\iota \co  H_1(M) \to    
 \Hom(H_1(M) ,\Qq/\Zz)$. Let $H_1(M)^*$ denote $ \Hom(H_1(M), \mathbb{C}^*). $ 
Let $\nu$ denote the map  $\Qq/\Zz \rightarrow \mathbb{C}^*$ that sends $\frac a b $ to $e^{\frac {2 \pi i a}{b}} .$
So we have an isomorphism $\jmath \co  H_1(M) \to H_1(M)^*$ given by $x \mapsto  \nu \circ \iota(x).$ Let $\beta \co   H_1(M)^* \times H_1(M)^* \rightarrow \Qq/\Zz $ be the dual form
 defined by    
 $ \beta(\jmath x, \jmath y)= -  l(x,y)$.    
   
\begin{definition}   
The form $\beta$ is metabolic with metabolizer $H$ if there exists a subgroup   
 $H$ of $H_1(M)^*$ such that $H^\bot = H$.   
\end{definition}    
   
\begin{lemma}  
\label{extend}  
 {\rm \cite{Gi1}}\qua  
 If $M$ bounds a  
spin  
 $4$-manifold $W$ then $\beta = \beta_1 \oplus \beta_2$   
 where $\beta_2$ is metabolic and $\beta_1$ has an even presentation with rank $\text{dim } H_2(W;\Qq)$  
 and 
signature Sign($W$).   
 Moreover, the set of characters that extend to $H_1(W)$  forms a metabolizer   
 for $\beta_2$.    
\label{linking}   
\end{lemma}   
  
\subsection{Link invariants}  \label{LI}
  
Let $L=L_1 \cup \dots \cup L_{\mu}$ be an oriented link in $S^3$.  
Let $N_2$ be the two-fold covering of $S^3$ branched along $L$ and $\beta_L$ be  
 the linking form on $H_1(N_2)^*$, see previous section.  
  
 We suppose that the Alexander polynomial of $L$ 
satisfies
 $$ \Delta_L(-1) \neq 0. $$  
Hence, $N_2$ is a rational homology sphere. Note that if $ \Delta_L(-1) \neq 1$,  
 then $H_1(N_2;\Zz)$ is non-trivial. 
  
\begin{definition}  
For all 
characters
 $\chi$ in  $H_1(N_2)^*$, the Casson-Gordon $\sigma$-invariant  
 of $L$ and the related nullity are (see Definition \ref{sigma}):  
  $$ \sigma(L,\chi) := \sigma(N_2,\chi), $$  
$$ \eta(L,\chi):= \eta(N_2,\chi). $$  
\label{cg}  
\end{definition}

\begin{remark}  
If $L$ is a knot, then Definition \ref{cg} coincides with $\sigma(L,\chi)$  
 defined in \cite[p.183]{CG1}.  
\end{remark}  
  
\section{Framed link descriptions}  
  
In this section, we study the Casson-Gordon $\tau$-invariants of the two-fold 
 cover $M_2$  
 of the manifold $M_0$  
 described below.
  
Let $S^3 - T(L)$ be the complement in $S^3$ of an open tubular neighborhood of $L$   
 in $S^3$ and  $P$ be a planar surface with $\mu$ boundary components.  
  
Let  $S$ be a Seifert surface for $L$ and $\gamma_i$ for $i=1, \dots, \mu$ 
be  
 the curves where $S$ intersects the  boundary of  $S^3 - T(L)$.   
 We define  
 $M_0$ as the result of gluing $P \times S^1$ to $S^3 - T(L)$, where $P \times {1}$   
 is glued along the curves $\gamma_i$. Let $*$ be a point in the boundary of $P$.  
  
A recipe for drawing a framed link description for $M_0$ is given in the proof of  
 Proposition \ref{M0}.  
   
\begin{proposition}   
$$H_1(M_0) \simeq \Zz \oplus \Zz^{\mu -1} \simeq   \langle  m \rangle \oplus \  \Zz^{\mu -1},$$   
 where $m$ denotes the class of $*  \times S^1$ in $P \times S^1$.     
\label{M0}  
\end{proposition}   
   
\begin{proof}   
Form a 4-manifold $X$ by gluing $P \times D^2$ to $D^4$ along $S^3$  
in such a way that the total   
 framing on $L$ agrees with the Seifert surface $S$. The boundary of this   
4-manifold is $M_0$.   
 We can get a surgery description of  $M_0$ in
the following way:  
  pick $\mu - 1$ paths of $S$ joining up the components of $L$ in a chain.   
 Deleting   
open neighborhoods of  
these paths 
in
  $S$ gives a Seifert surface for a knot $L^\prime$ obtained   
 by doing a fusion of $L$ along bands that are neighborhoods of the original
paths.
 Put a circle with a dot around each of these bands (representing   
  a
4-dimensional $1$-handle in Kirby's \cite{K} notation), and the framing zero on $L^\prime.$  
This describes a handlebody decomposition of $X.$  
  
 One can then get    
 a standard framed link description of $M_0$ by replacing   
 the circle with dots with
 unknots $T_1,\dots,T_{\mu-1}$ framed zero.  
 This changes the $4$-manifold but not the boundary. Note also that   
 $lk(T_i,T_j)=0$ and $lk(T_i,L^\prime)=0$ for all $i=1, \dots, \mu-1$.   
   Hence $H_1(M_0) \simeq \Zz^\mu$ and $m$ represents one of the generators.   
\end{proof}   
   
We now consider an infinite cyclic covering $M_\infty$ of   
 $M_0$, defined by a character   
 $H_1(M_0) \to C_\infty= \langle  t   \rangle $ that sends $m$ to $t$ and the other  
generators  
to zero.  
 Let us denote  by $M_2$  the intermediate two-fold 
covering   
 obtained by composing this character with the quotient map $C_\infty \to C_2$ sending $t$ to $-1$.   
Let $m_2$ denote the loop  in $M_2$  given by the inverse image of $m$. A recipe for drawing a framed link description for $M_2$ is given in the proof of Remark \ref{surgM2}.   

\begin{proposition}  \label{M2}   There is an isomorphism between $H_1(N_2)$ and the torsion subgroup of  
$H_1(M_2)$, which only depends on $L.$ Moreover 
$$ H_1(M_2) \simeq H_1(N_2) \oplus \Zz^\mu    
 \simeq H_1(N_2) \oplus   \langle  m_2 \rangle \oplus \Zz^{\mu-1}. $$  
\end{proposition}   
   
\begin{proof}  

Let  $R$ be the result of gluing $P \times D^2$  to $S^3 \times I$ along  $L\times {1} \subset S^3 \times {1}$ using the framing given by the Seifert surface. Thus $R$ is the result of adding
 $\mu-1$ 1-handles to  $S^3 \times I$  and then one 2-handle along $L'$,  as  in the proof above.
Then $X$ in the proof above can be obtained by gluing $D^4$ to $R$ along $S^3 \times {0}.$ 
 Since $D^2$ is the double branched cover of itself along the origin, $P \times D^2$
 is the double branched cover of itself along $P \times {0}$.
 Let $R_2$ 
denote
 the double branched cover of
 $R$ that is obtained by gluing $P \times D^2$  to $N_2 \times I$ along  a neighborhood 
of
 the lift of $L\times {1} \subset S^3 \times {1}.$  We have that $\partial R_2= -N_2 \sqcup M_2$,  
 where $R_2$ is the result of adding $\mu-1$ 1-handles to
$N_2 \times I$ and then one 2-handle along  the lift $L'.$  Moreover this lift of $L'$  is
 null-homologous in $N_2.$  It follows that $H_1(R_2)$
is isomorphic to  $H_1(N_2) \oplus \Zz^{\mu -1},$  with the inclusion of $N_2$ into $R_2$  inducing an  isomorphism
 $i_N$ of 
$H_1(N_2)$ to  the torsion subgroup of $H_1(R_2).$  Turning this handle decomposition upside down we have that
$R_2$ is the result of adding to
$M_2 \times I$  one 2-handle along  a neighborhood of $m_2$ and then $\mu-1$ 3-handles. 
  It follows that $H_1(R_2) \oplus \Zz = H_1(R_2) \oplus    \langle  m_2 \rangle $
is isomorphic to  $H_1(M_2)$  with the inclusion of $M_2$ in $R_2 $ inducing  an isomorphism  $i_M$ of  the torsion subgroup $H_1(M_2)$ to  the torsion subgroup of $H_1(R_2).$  Thus  $(i_M)^{-1} \circ  i_N$ is an isomorphism
from $H_1(N_2)$ to the torsion subgroup of  $H_1(M_2)$ 
and  this isomorphism  is constructed
without 
any arbitrary choices.\end{proof}  

\begin{remark} \label{surgM2}
We could have proved Proposition \ref{M0} in a similar way to the proof  of   Proposition  \ref{M2}. 
We could have 
also
proved   Proposition  \ref{M2} (except for the isomorphism only depending on $L$) in a similar way to the proof  of Proposition  \ref{M0} as follows.
 We can  find a surgery description of $M_2$ from    
a surgery description of $N_2$.   
 The procedure of how to visualize a lift of $L$ and the surface $S$   
 in $N_2$ is given in \cite{AK}. One considers the lifts of  the paths chosen in 
 the proof of Proposition \ref{M0}, 
on the lift of $S.$ One 
then
fuses the components of the lift of $L$ along these paths,
 obtaining a lift of $L'.$
The surgery description of $M_2$    
 is obtained by adding 
 to the surgery description   
 of $N_2$
the lift of $L'$ with zero framing  together with $\mu-1$ more unknotted zero-framed components 
encircling
 each fusion.  The linking matrix of this   
 link is a direct sum of that of $N_2$ and  
a $\mu \times \mu$ zero matrix. \end{remark}
\eject

 Let $i_T$ denote the inclusion of  the torsion subgroup of  $H_1(M_2)$ into  $H_1(M_2),$   and let
 $\psi \co  H_1(N_2) \rightarrow H_1(M_2)$ denote the monomorphism given by 
  $i_T \circ (i_M)^{-1} \circ  i_N.$
\begin{theorem}  
  
Let $\chi^+ \co H_1(M_2) \to \Cc_q \oplus C_\infty.$  Let $\chi \co H_1(N_2) \to \Cc_q $  
 be $\chi^+ \circ \psi$ composed with the projection to  $\Cc_q . $  We have that:  
$$ | \sigma_1(\tau(M_2,\chi^+)) -\sigma(L,\chi)| \leq  
 \eta(L,\chi) + \mu.$$  
\label{est}  
\end{theorem}  
  
\begin{remark}  
If $L$ is a knot, then $\tau(M_2,\chi^+)$ coincides with 
 $\tau(L,\chi)$
 defined  
 in \cite[p.189]{CG1}.  
\end{remark}  

 \proof[Proof of Theorem \ref{est}]   
We use the surgery description of $M_2$ given in Remark 
 \ref{surgM2}.  Let $P$ be given by the surgery description  
of $M_2$ but with the component corresponding to $L^\prime$ deleted.  
 Hence, $$ P = N_2 \sharp_{(\mu - 1)} S^1 \times S^2. $$  
$\chi+$ 
 induces 
some
character  $\chi^\prime$ on $H_1(P)$.

 According to Section \ref{cgt}, 
 we let
 $\overline{\chi} \in H^1(M_2;\Cc_q)$ and
 $\overline{\chi}^\prime \in H^1(P;\Cc_q)$ 
denote
 the
 characters $\chi^+$ and $\chi^\prime$  followed by the projection $\Cc_q \oplus C_\infty \to \Cc_q$.
Using Propositions \ref{add} and \ref{s1s2}, one has that 
$$ \sigma(P,\overline{\chi}^\prime) = \sigma(L,\chi) \text{ and }  \eta(P,\overline{\chi}^\prime)=\eta(L,\chi) + \mu - 1. $$ 
 Moreover, since $M_2$ is obtained by surgery on $L^\prime$ in $P$, it     
 follows from \cite[Proposition (3.3)]{Gi3} that  
  $$ |\sigma(P,\overline{\chi}^\prime) - \sigma(M_2,\overline{\chi}) |+|  \eta(M_2,\overline{\chi}) - 
 \eta(P,\overline{\chi}^\prime) | \leq  1 \text{ or }$$  
$$ |\sigma(L,\chi) - \sigma(M_2,\overline{\chi}) |+|  \eta(M_2,\overline{\chi}) -  \eta(L,\chi) - \mu +1| \leq  1.$$  
Thus  
$$ |\sigma(L,\chi) - \sigma(M_2,\overline{\chi}) | \leq   \eta(L,\chi) + \mu - \eta(M_2,\overline{\chi}) .$$  
Finally, one gets, by Theorem \ref{firstest},  
$$ | \sigma_1(\tau(M_2,\chi^+)) -\sigma(L,\chi)| \leq | \sigma_1(\tau(M_2,\chi^+)) - \sigma(M_2,\overline{\chi}) |  
 + | \sigma(M_2,\overline{\chi}) - \sigma(L,\chi)|  $$  
$$ \leq \eta(M_2,\overline{\chi})+ \eta(L,\chi) + \mu - \eta(M_2,\overline{\chi}) =  \eta(L,\chi) + \mu.\eqno{\qed} $$  

\eject

\section{The slice genus of links}  
 
See Section \ref{LI} for notations.

\begin {theorem} \label{main} 
Suppose $L$ is the boundary of a connected oriented 
properly embedded
surface $F$ of  genus $g$ in $B^4,$  
and that $\Delta_L(-1) \neq 0$.  Then,  
$\beta_L$ can be 
written as a direct sum $\beta_1 \oplus \beta_2$ such that the following two conditions hold:   
  
{\rm 1)}\qua $\beta_1$   has an even presentation of    
rank $2g + \mu -1$  and signature $\sigma_L(-1)$, and $\beta_2$ is metabolic.   

{\rm 2)}\qua There is a metabolizer for $\beta_2$  such that for all characters $\chi$ of prime power order in this metabolizer,   
  $$ | \sigma(L,\chi)+\sigma_L(-1) | \leq  \eta(L,\chi) + 4g +3 \mu -2 .$$  
\end{theorem}  
  
\proof    
We let $b_i(X)$ denote the ith Betti number of a space $X$.  
 We have $b_1(F) = 2g + \mu -1.$

Let $W^\prime_0$, with boundary $M^\prime_0$, be the complement of an open tubular neighborhood of   
 $F$  
 in $B^4$. 
By the Thom isomorphism, excision, and the long exact sequence of the pair 
$(B^4,W^\prime_0),$ $W^\prime_0$
has the homology of $S^1$ wedge 
$b_1(F)$ 2-spheres.
Let $W^\prime_2$ with boundary $M^\prime_2$ be the two-fold  
covering of $W^\prime_0$. Note that    
if $F$ is planar,  $M^\prime_0= M_0, $  and $M_2^\prime=M_2$  
 (see Section 3).

  Let $V_2$ be the two-fold covering of $B^4$ with branched set $F$.   
Note that  
$V_2$ is spin as  $w_2(V_2)$ is the pull-up 
 of a class in $H^2(B^4,\Zz_2)$,  by  \cite[Theorem 7]{Gi5}, for instance.  
The boundary of $V_2$ is $N_2$.
As in  \cite{Gi1}, one calculates that
 $b_2(V_2)=2g + \mu -1$.
One has  $\Sign(V_2)=  \sigma_{L}(-1)$ by \cite{V}.

By Lemma {\ref{extend}}, $\beta_L$ can be written as a direct sum $\beta_1 \oplus \beta_2$ as in condition $1)$ above, such that the characters on $H_1(N_2)$ that extend   
 to $H_1(V_2)$ form a metabolizer $H$ for $\beta_2$. We now suppose $\chi \in H$ and show that Condition $2)$ holds for $\chi.$  
  
We also  let  $\chi$ 
 denote
 an extension of $\chi$ to $H_1(V_2)$  
 with image some cyclic group $\Cc_q$ where $q$ is a power of a prime  
integer  
 (possibly   
 larger than those corresponding to the character on $H_1(N_2)$). 
 Of course $\chi \in H^1(V_2,\Cc_q) $ restricted to $W_2^\prime$ extends  
$\chi $ restricted to $M_2^\prime$.  
 We simply  denote all these restrictions by $\chi$.

 Let $W^\prime_\infty$ denote the infinite cyclic cover of  
$W^\prime_0$. Note that $W^\prime_2$ is a quotient of this covering space.
 $\chi$ induces a $\Cc_q$-covering of $V_2$ and thus of $W_2^\prime$.   
 If we pull the $\Cc_q$-covering
 of $W_2^\prime$ up to $W_\infty^\prime$, we obtain   
 $\widetilde{W}_\infty^\prime$,  a $\Cc_q \times C_\infty$-covering 
 of $W_2^\prime$.   
 If we identify properly $F \times S^1$ in $M_2^\prime,$ this covering restricted to   
 $F \times S^1$ is given by a character $H_1(F \times S^1) \simeq H_1(F) \oplus   
 H_1(S^1) \to \Cc_q \times C_\infty$ that maps $H_1(F)$ to zero in   
 $C_\infty$, $H_1(S^1)$ to zero in $\Cc_q$ and isomorphically onto   
 $C_\infty$.  
For this note: since  
 $\text{Hom}(H_1(F),\Zz)=H^1(F)= [F,S^1]$,  
  we may define diffeomorphisms  
of $F \times S^1$ that induce the identity on the second factor  
of $H_1(F \times S^1) \approx H_1(F) \oplus \Zz,$ and send 
$(x,0)\in H_1(F) \oplus \Zz,$ to $(x,f(x))\in H_1(F) \oplus \Zz,$  for any 
  $f \in \text{Hom}(H_1(F),\Zz).$
  
As in \cite{Gi1},  
 choose inductively a collection of $g$ disjoint curves in 
the kernel of  
$\chi$   
 that form a metabolizer for the intersection form  on 
$H_1(F)/ H_1( \partial F)$.
  By taking a tubular neighborhood of these curves   
 in $F$, we obtain a collection of $S^1 \times I$ embedded in $F$. Using these   
 embeddings we can attach round 2-handles   
 $(B^2 \times I) \times S^1$ along $(S^1 \times I) \times S^1$ to the trivial   
 cobordism  
$M_2^\prime \times I$ and obtain a cobordism $\Omega$ between  $M_2$ and $M_2^\prime$.   
   
 Let $U= W_2^\prime \cup_{M_2^\prime} \Omega$
 with  
 boundary $M_2$.  
The $\Cc_q \times C_\infty$-covering of $W^\prime_2$  extends uniquely to $U$.  
 Note that
 $\Omega$ may also be viewed as the result of attaching round 1-handles to $M_2 \times I.$   

 As in \cite{Gi1},  ${\Sign}(W^\prime_2)= {\Sign}(V_2)$.   
Since the intersection form on  $\Omega$ is zero,  we get  
  ${\Sign}(U)={\Sign}(W_2^\prime)= {\Sign}(V_2)=\sigma_L(-1)$.  
 The $\Cc_q \times C_\infty$-covering of $\Omega$, restricted to each round $2$-handle   
 is $q$ copies of $B^2 \times I \times \Rr$ attached to the trivial cobordism $\widetilde{M}_\infty^\prime \times I$ 
along $q$ copies of $S^1 \times I \times \Rr$.   
 Using a Mayer-Vietoris sequence, one sees that
 the inclusion induces an isomorphism (which
 preserves    
 the Hermitian form) 
 $$ H_2^t(U;\Qq(\Cc_q)(t)) \simeq H_2^t(W_2^\prime;\Qq(\Cc_q)(t)). $$   
 Thus, if $w(W_2^\prime)$ denotes the image of the intersection form on $H_2^t(W_2^\prime;\Qq(\Cc_q)(t))$   
 in $\mathcal{W}$ $(\Qq(\Cc_q)(t))$, we get $\sigma_1(\tau(M_2, \chi^+)) = \sigma_1(w(W_2^\prime))-\sigma_{L}(-1)$.   
     
 If $q$ is a prime power, we may apply Lemma 2 of \cite{Gi1} and  
conclude  
that   
 $H_i(\widetilde{W}_\infty^\prime;\Qq)$
 is finite dimensional for all $i \neq 2$.   
  Thus,  $H_i^t(W_2^\prime; \Qq(\Cc_q)(t))$ is zero for all $i \neq 2$. Since   
 the Euler characteristic of 
$W_2^\prime$
 with 
 coefficients  in $\Qq(\Cc_q)(t)$ coincides   
 with those with 
 coefficients in $\Qq$, we get $\text{dim} \ H_2^t(W_2^\prime;\Qq(\Cc_q)(t)) = \chi(W_2^\prime)= $
$ 2\chi(W_0^\prime)= 2(1-\chi(F))= 2 b_1(F)$. 
Thus 
$| \sigma_1 (\tau (M_2,\chi^+)+\sigma_L(-1) | \leq  2 b_1(F)$.  
Hence,  
$$ | \sigma(L,\chi) +  \sigma_L(-1) | \leq | \sigma(L,\chi) - \sigma_1 (\tau (M_2,\chi)^+)| + | \sigma_1 (\tau (M_2,\chi)^+)  
 + \sigma_L(-1) | $$  
$$ \leq \eta(L,\chi) + \mu + 2 ( 2g + \mu - 1) = \eta(L,\chi) + 4g + 3 \mu - 2 \text{ by Theorem \ref{est}. }\eqno{\qed} $$

\section{Examples}  

Let $L=L_1 \cup L_2$ be the link with two components of Figure 1   
 and $S$ be the Seifert surface of $L$ given by the picture.    
 The 
 squares with $K$ denote two parallel copies with linking number $0$ of 
an arc tied in
the knot $K$.
Note that $L$ is actually a family of examples. Specific links are determined by the choice of two parameters: a knot $K$ and a positive integer $h.$
 Since $S$ has genus $h$,  the slice genus 
of $L$
is at most $h$.  

\begin{figure}[htbp]   
\begin{center}   
\includegraphics[width=8.4cm,height=5.4cm]{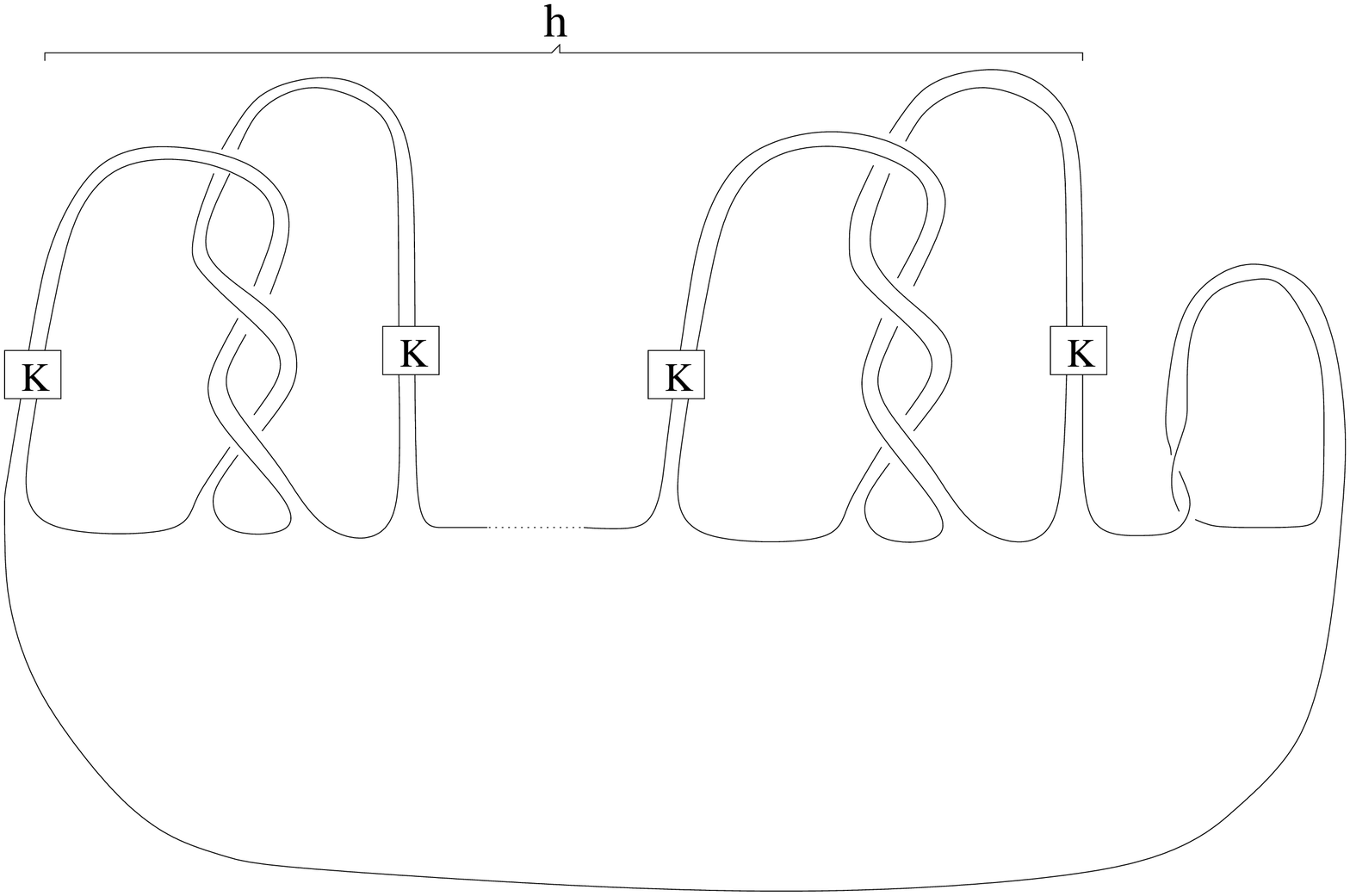}   
\end{center}   
\caption{The link $L$}   
\end{figure}

One calculates that $\sigma_{L}(\lambda) =1$, and  $n_{L}(\lambda) = 0$ for all $\lambda$. 
 Thus, the Murasugi-Tristram  inequality
says nothing about the slice genus of $L$. 
In fact, if $K$ is a slice knot, then one can surger this surface to obtain a smooth cylinder
 in the 4-ball with boundary $L$. Thus there can be no arguments based solely on a Seifert pairing for $L$
that would imply that the slice genus is non-zero.   
\begin {theorem} \label{example} If  
 $\sigma_K(e^{2 i \pi / 3}) \ge 2h$ or  $\sigma_K(e^{2 i \pi / 3}) \le -2h -2,$  then $L$ has slice genus $h$.   
  \end{theorem}

\begin{proof}  
  Using \cite{AK}, a surgery presentation of $N_2$ as surgery on a framed link of $2h+1$ components can be obtained from   
 the surface $S$ (see Figure 2).

 Let $Q$ be the $3$-manifold obtained from the link pictured in
 Figure 2.  
Here $K'$ denotes $K$ with the string orientation reversed.
Since $RP(3)$ is  obtained by surgery on the unknot framed $2$, we get:
$$ N_2 = RP(3) \#_{h} Q. $$
     
\begin{figure}[htbp]   
\begin{center}   
\includegraphics[width=6cm,height=3.6cm]{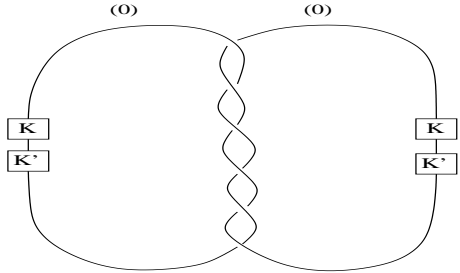}   
\end{center}   
\caption{Surgery presentation of $Q$}   
\end{figure}   
   
The linking matrix of 
 the framed link of the surgery presentation of $N_2$
 is  

$\Lambda = [2] \bigoplus \oplus^h \begin{bmatrix} 0 & 3 \\   
3 & 0  \end{bmatrix}$.
  $\Lambda$  
 is a presentation matrix of $(H_1(N_2)^*,\beta_{L})$; we obtain  
   $$ H_1(N_2)^* \simeq \Zz_2 \bigoplus \oplus^{2h} \Zz_3 $$
 and $\beta_{L}$ is given by the following matrix, with entries in $\Qq / \Zz$:  $$[1/2] \bigoplus \oplus^h \begin{bmatrix} 0 & 1/3 \\   
1/3 & 0  \end{bmatrix}.$$  
  By Theorem \ref{main}, if ${L}$ bounds 
a
surface of genus $h-1$ in $B^4$, then $\beta_{L}$ must   
be decomposed as $\beta_1 \oplus \beta_2$ where:   
   
1)\qua $\beta_1$ has an even presentation matrix of rank 
$2h -1$, 
 and signature $1$  (all we really 
need
 here is that it has a rank
$2h-1$
 presentation.)  
   
2)\qua $\beta_2$ is metabolic and for all characters $\chi$ of prime power order in  
some metabolizer of $\beta_2$,  
 the following inequality holds:  
$$ | \sigma({L},\chi) + 1 | - \eta({L},\chi) \leq 4h   .
  \leqno{(*)}$$
 As 
$\Zz_2 \bigoplus \oplus^{2h} \Zz_3 $
 does not have a rank 
 $2h -1$
 presentation, $\beta_2$ is non-trivial. 
As metabolic forms are defined on groups whose cardinality is a square,
 $\beta_2$ is defined on a group with no $2$-torsion.
Thus the metabolizer contains a non-trivial 
 character of order three satisfying  $\beta_{L}(\chi,\chi) =0.$ 

The first homology of $Q$ is  $\Zz_3 \oplus  \Zz_3$, generated by, say, $m_1$ and $m_2$, positive meridians of these components.
Each of these components is oriented counterclockwise.
We first work out $\sigma(Q,\chi)$ and  $\eta(Q,\chi)$ for characters of order three. Let $\chi_{(a_1,a_2)}$ denote the character on
$H_1(Q)$ sending 
$m_j$
 to  $e^{\frac{2i \pi a_j}{3}}$,  where the $a_j$
 take 
the
values zero and $\pm1.$

We use Proposition \ref{surgeryformula} to compute  $\sigma({Q},\chi_{(1,0)})$ and  $\eta({Q},\chi_{(1,0)})$ assuming that $K$ is trivial. For this, one may adapt the trick illustrated on a link with  
$2$  twists between the components \cite[Fig (3.3), Remark (3.65b)]{Gi2}. 
 In the case $K$ is the unknot, we obtain 
$$ \sigma({Q},\chi_{(1,0)}) = 1 \quad \text{ and } \quad  \eta({Q},\chi_{(1,0)}) = 0.$$
It is not difficult to see that inserting 
the knots of the type $K$ changes the result as follows 
(note that $K$ and $K^\prime$ have the same 
Tristram-Levine  signatures):
$$\sigma(Q,\chi_{(1,0)}) = 1 + 2 \sigma_K(e^{2 \pi i / 3}) \quad \text{ and } \quad  \eta(Q,\chi_{(1,0)}) = 0.$$  
These same values hold for the characters   
$\chi_{(-1,0)}$
 and $\chi_{(0,\pm 1)}$ by symmetry.

Using Proposition \ref{surgeryformula} 
$$\sigma(Q,\pm \chi_{(1,1)}) = -1 -24/9 + 4 \sigma_K(e^{2 \pi i / 3}) , \quad  \quad  \eta(Q,\pm \chi_{(1,1)}) = 0$$
$$\sigma(Q,\pm \chi_{(1,-1)}) = 4 + 24/9 + 4 \sigma_K(e^{2  \pi  i/ 3}) \quad \text{ and } \quad  \eta(Q,\pm \chi_{(1,-1)}) = 1.$$  
One also has 
$$\sigma(Q,\chi_{(0,0)}) = 0 \quad \text{ and } \quad  \eta(Q,\chi_{(0,0)}) = 0.$$
Any order three character on $N_2$ that is self annihilating under the linking form is given as the sum of
the trivial character on $RP(3)$ and characters of type $\chi_{(0,0)}$,  $\chi_{(\pm 1,0)}$ and $\chi_{(0,\pm 1)}$  on $Q$ and characters of type $\pm \chi_{(1,1)} + \pm \chi_{(1,-1)}$ on $Q \# Q$.
Using Proposition \ref{add}, one can calculate $\sigma(L,\chi)$ and $\eta(L,\chi)$ for all these characters $\chi$. It is now a trivial matter to check that for every 
 non-trivial
 character 
with $\beta (\chi,\chi)=0$, 
the inequality (*) is not   
satisfied.
\end{proof}

   \Addresses  
    \end{document}